\documentclass[12pt]{article}
\usepackage[a4paper]{geometry}
\usepackage{amssymb,amsbsy}
\usepackage{times}

\usepackage{datetime}

\newcommand{\cL}{{\cal L}}

\newcommand{\RR}{\mathbb{R}}
\newcommand{\NN}{\mathbb{N}}

\newcommand{\Span}{{\mathop{\rm span\,}\nolimits}}

\newcommand{\eop}{\hfill$\Box$}

\newtheorem{theorem}{Theorem}[section]
\newtheorem{proposition}[theorem]{Proposition}
\newtheorem{lemma}[theorem]{Lemma}
\newtheorem{corollary}[theorem]{Corollary}

\newtheorem{remark}[theorem]{Remark}
\newtheorem{definition}[theorem]{Definition}

\newenvironment{pf}%
{\par\noindent\textbf{Proof:}~}%
{\eop\par\smallskip\par\noindent}

{\par\noindent\textbf{Proof of #1:}~}%
{\eop\par\smallskip\par\noindent}

\begin{document}

\title{Observations on interpolation by total degree polynomials in
  two variables} 
\author{Jes\'us Carnicer\thanks{Departamento de Matematica
    Aplicada/IUMA, Universidad de Zaragoza, Pedro Cerbuna 12, 50009
    Zaragoza, Spain. \texttt{carnicer@unizar.es}. Partially supported
    by MTM2015-65433-P (MINECO/FEDER) Spanish Research Grant, by Gobierno de Arag\'on
    and European Social Fund.}
  \and Tomas Sauer\thanks{Universit\"at Passau, Lehrstuhl f\"ur Mathematik mit
    Schwerpunkt Digitale Bildverarbeitung \& FORWISS, Innstr. 43,
    D-94032 Passau, Germany. \texttt{Tomas.Sauer@uni-passau.de}}}
\maketitle

\begin{abstract}
  In contrast to the univariate case, interpolation with polynomials
  of a given maximal total degree is not always possible even if the
  number of interpolation points and the space dimension coincide. Due
  to that, numerous constructions for interpolation sets have been
  devised, the most popular ones being based on intersections of
  lines. In this paper, we study algebraic properties of some such
  interpolation configurations, namely the approaches by
  Radon-Berzolari and Chung-Yao. By means of proper H-bases for the
  vanishing ideal of the configuration, we derive properties of the
  matrix of first syzygies of this idealwhich allow us to draw
  conclusions on the geometry of the point configuration.
\end{abstract}

\section{Introduction}
\label{sec:Intro}

Interpolation of data, especially by polynomials, is a classical
issue, not only in one
variable but also in several variables,
cf. \cite{GascaSauer00b,Steffensen27}. While in one variable the
interpolation polynomial can be easily expressed in a closed form,
this is not the case in two and more variables where the geometry becomes
significantly more intricate, especially due to the fact that there are
no multivariate Haar spaces, see \cite[Chapter 2, Section 4]{Lorentz66}.

To overcome some of these problems, techniques were developed to
construct point sets that allow for unique interpolation from a
given subspace of polynomials, typically the vector space of all
polynomials of total degree not greater than a given number.
The Radon--Berzolari construction \cite{berzolari14:_sulla,Radon48}
generates a set of interpolation points by decomposing a bivariate
problem of 
degree, say $n$, into two simpler subproblems, one being univariate
and one of degree $n-1$: the construction
extends a set of interpolation points for degree $n-1$ by choosing $n+1$
additional points on a line. As long as this line does not contain any
of the low degree interpolation points, this gives a valid set of
interpolation points of degree $n$.

Later, Chung and Yao \cite{ChungYao77} presented a \emph{geometric
  characterization} of certain interpolation sets that extend
properties of the univariate case to several variables. In particular,
they defined a class of interpolation sets which are nowadays known as \emph{GC
configurations} \cite{CarnicerGasca04} or \emph{GC sets}
\cite{boor07:_multiv,carnicer06:_geomet}. Gasca and Maeztu
\cite{GascaMaeztu82} conjectured that any \emph{bivariate} GC set is
the result of
a Radon--Berzolari construction, that is, any such set possesses a 
\emph{maximal line} that contains $n+1$ of the interpolation points.
This conjecture, based on the simple observation in the cases $n=1,2$, has
so far only been proven for degree up to $5$, see
\cite{hakopian14:_gasca_maezt}.

Recently, Hal Schenck pointed out the striking
connections between
interpolation sets and the generating matrix for the first syzygy
module of the associated zero dimensional radical ideal. This
approach which is based on sophisticated
concepts from Algebraic Geometry, see
\cite{eisenbud05:_geomet_syzyg}, can be found in \cite{Schenck16P}.
One main point there is that it is possible
to characterize the existence of a maximal line in an 
interpolation set by looking at the syzygy matrix of the respective
ideals. The purpose of this paper is to give
an elementary and direct approach to these ideas in which the
Berzolari--Radon construction plays a significant role.

\section{Interpolation of total degree and  ideals}

Let  $\Pi= \RR[x]$, $x=(x_1,x_2)$, be the set of bivariate polynomials
and  $\Pi_n$ be the set of bivariate polynomials of  total degree less
than or equal to $n$ whose dimension is $ \dim \Pi_n =(n+1)(n+2)/2$.
Given a  set $Y \subseteq \RR^2$, the evaluation map is defined
as 
$$
 p \in \Pi \mapsto p(Y):=(p(y))_{y \in Y} \in \RR^Y.
$$
The kernel of the evaluation map is the \emph{ideal} 
$$
I(Y):=\{ p \in \Pi :p(Y) =0\}
$$
of all bivariate polynomials vanishing at the set $Y$. 

Given  a \emph{finite} set of  nodes $Y\subset \RR^2$ and $f\in\RR^Y$,
we can formulate the \emph{Lagrange interpolation problem} on $Y$ in
$\Pi_n$: find $p\in \Pi_n$ such that  $p(Y) = f$. 

\begin{definition}\label{D:PoisedDef}
A set $Y \subset \RR^2$ is  $\Pi_n$-\emph{independent} if the
Lagrange interpolation problem on $Y$ has a solution in $\Pi_n$, maybe not
unique.  A set  $Y \subset \RR^2$ is $\Pi_n$-\emph{poised} if the
interpolation problem on  $Y$ in $\Pi_n$ is unisolvent, that is, the
interpolation problem on $Y$ has always a \emph{unique} solution in $\Pi_n$.
\end{definition}

The  evaluation map is a surjective linear map for any finite $Y$
since for each $f \in \RR^Y$ the polynomial
$$
p(x) = \sum_{y \in Y} f(y) \prod_{t \in Y\setminus\{y\}} {(x_1-t_1)^2
  + (x_2-t_2)^2 \over (y_1-t_1)^2 + (y_2-t_2)^2}
$$
satisfies $p(Y)= f$.  Hence, $\RR^Y \cong \Pi/I(Y)$.  

The set $Y$ is $\Pi_n$-poised if and only if the restriction of the
evaluation map to $\Pi_n$ is bijective.  So, if $Y$ is $\Pi_n$-poised then
$\dim \Pi_n = \dim \RR^Y= \# Y$.    Uniqueness of the solution of the
Lagrange interpolation problem implies that $\Pi_n \cap I(Y)=0$ and
existence of a solution implies that $\Pi= \Pi_n + I(Y)$.  Therefore
$Y$ is $\Pi_n$-poised if and only if $\Pi_n \oplus  I (Y)= \Pi$.

If $Y$ is a $\Pi_n$-poised set, then for each $y \in Y$ there exists a
unique \emph{fundamental polynomial}, also called \emph{Lagrange polynomial},
$\ell_{y,Y}$ in $\Pi_n$ such that
$$
\ell_{y, Y}(y') = \delta_{yy'}, \quad y' \in Y,
$$
where $\delta_{yy'}$ is Kronecker symbol.  For a  $\Pi_n$-poised set
$Y$ the Lagrange interpolation operator $L_{Y}$ associates to each
function from $\RR^2$ to $\RR$ its  polynomial interpolant in
$\Pi_n$. Considered as an operator from $\Pi$ to $\Pi$,
the Lagrange interpolation operator is a projection that can be
expressed in terms of the fundamental polynomials by means of the
Lagrange formula 
$$
L_{Y} [f]=\sum_{y \in Y} f(y) \ell_{y,Y}.
$$
The \emph{error operator} 
$$
E_{Y}[f] := f - L_{Y}[f]
$$ is another linear
projection whose image is the ideal $I(Y)$, and the two projectors are
complementary.  Because of that the Lagrange interpolation operator is
called an \emph{ideal projector}, that is, a projector whose kernel is
an ideal \cite{boor05:_ideal}.

From the definition it follows that $Y$ is a $\Pi_n$-independent set
if and only if  for each $y \in Y$ there exists a fundamental
polynomial (maybe not unique) in $\Pi_n$ vanishing at $Y \setminus
\{y\}$ and with value 1 at  $y$. Observe that a subset of a
$\Pi_n$-independent  set is also $\Pi_n$-independent. 
For each $\Pi_n$-independent set we have that  $\# Y \le \dim \Pi_n$
and, if equality holds, then $Y$ is $\Pi_n$-poised. 

\begin{definition}
  A subset $H$ of an ideal $I$ is called an
  \emph{H-basis} for $I$ if any $f \in I$ can be written as
  $$
  f= \sum_{h \in H} g_h h, \qquad g_h \in \Pi_{\deg f- \deg h}, \quad h \in H.
  $$  
\end{definition}

\begin{lemma}\label{L:2.2}
  Let $I$ be an ideal of $\Pi$ such that $\Pi_{n}\cap I = 0$ and
  $\dim ( \Pi_{n+1} \cap I)=n+2$. 
  Then $h_0, \ldots, h_{n+1} $ is a basis of $\Pi_{n+1} \cap I $ if
  and only if  it is an H-basis of $I$. Moreover, $\Pi=\Pi_n\oplus
  I$.
\end{lemma}

\begin{pf}
  Let $h_0, \ldots, h_{n+1}$ be a basis of $ \Pi_{n+1} \cap I$. Since
  $\Pi_{n}\cap I = 0$, we have that $\Pi_{n+1} = \Pi_n \oplus
  (\Pi_{n+1} \cap I)$ and for each $\alpha$ with $\vert
  \alpha \vert = n+1$ there exists  $g_{\alpha}$ in  the vector space
  $\Pi_{n+1} \cap I$  such that the $g_{\alpha}(x)-x^{\alpha} \in \Pi_n$. 
  Since  $x^\alpha$, $\vert \alpha \vert= n+1$, are linearly
  independent, the polynomials $g_{\alpha}(x)$, $\vert
  \alpha\vert=n+1$, are also linearly independent and, since $\dim
  (\Pi_{n+1} \cap I) = n+2$, they are a basis of
  $\Pi_{n+1} \cap I$. 
  Now take any polynomial $p \in \Pi_m$, $m \ge n+1$.
  The homogeneous leading form of $p$ can be expressed as
  $\sum_{\vert \alpha \vert  =n+1} c_\alpha (x)x^\alpha$ with
  $c_{\alpha}(x) \in \Pi_{m-n-1}$.  Then  
  $$
  p-\sum_{\vert \alpha \vert  =n+1} c_\alpha g_\alpha \in  \Pi_{m-1}.
  $$
  An inductive argument shows that this reduction process allows us
  find polynomials $b_\alpha \in \Pi_{m-n-1}$, $|\alpha| = n+1$, such that 
  $$
  r:= p-\sum_{\vert \alpha \vert  =n+1} b_\alpha g_\alpha \in  \Pi_n.
  $$ 
  Hence it follows that $\Pi=\Pi_n\oplus I$.  If $p \in I$ then $r \in
  \Pi_n \cap I=0$, hence $p$ can be expressed 
  as a combination of the polynomials $g_\alpha$ with polynomial
  coefficients
  $b_\alpha\in \Pi_{m-n-1}$, $\vert \alpha \vert=n+1$.
  So,  $(g_\alpha:\vert \alpha \vert= n+1)$ is an $H$-basis. Since
  each $g_\alpha \in \Pi_{n+1} \cap I$ can be expressed as a linear
  combination of the basis $(h_0, \ldots, h_n)$ with constant
  coefficients and vice versa, it follows that $(h_0, \ldots, h_n)$ is
  also an H-basis.

  Conversely, if $h_0, \ldots, h_{n+1}$ is an
  H-basis and $p \in \Pi_{n+1} \cap I$, then $p$ can be written as a
  combination of the $h_j$, $j=0, \ldots, n+1$, where the coefficients
  are polynomials of degree $0$.  Hence $\Span\{h_0, \ldots, h_n\} =
  \Pi_{n+1} \cap I$ and, since $\dim (\Pi_{n+1} \cap I) = n+2$, it
  follows that $h_0, \ldots, h_n$ is a basis of the vector space $
  \Pi_{n+1} \cap I$.
\end{pf}

\noindent
The statement of the preceding lemma can be rephrased as that for any
ideal $I$ of $\Pi$ the conditions
$\Pi_{n}\cap I = 0$ and $\dim (\Pi_{n+1} \cap I)=n+2$ are
equivalent to $\Pi=\Pi_n\oplus I$.

\begin{corollary}\label{C:HBasisContained}
  Let $I$ be an ideal of $\Pi$ such that $\Pi_{n}\cap I = 0$ and
  $\dim ( \Pi_{n+1} \cap I)=n+2$ and let $P$ be a finite spanning subset of
  $\Pi_{n+1}\cap I$.  Then $P$ is an H-basis of $I$.
\end{corollary}

\section{Some remarks on the Berzolari-Radon construction}

For each set $K$ of bivariate polynomials, the associated algebraic
variety  is defined as $V(K):= \{x \in \RR^2 : k(x) = 0, \, k
\in K \}$.  If $K$ consists of a single polynomial $k$, we shall
denote by $K$ the algebraic curve with equation $k(x) = 0$ for the
sake of brevity, instead of using the  notation $V(K)$.  A line is the
set of zeros of a bivariate polynomial of first degree. 

Berzolari \cite{berzolari14:_sulla} and Radon \cite{Radon48} proposed
the  construction of a  $\Pi_{n+1}$-poised set $Y_{n+1}$ starting from
a $\Pi_{n}$-poised set  $Y_n$ by adding $n+2$ nodes lying on a line
$K$ that does not contain any node in $Y_n$. 
In (3) Fact of \cite{boor07:_multiv}, we can find a proof of this result  
and a relation between the fundamental polynomials of both sets.
We provide a restatement
of these results,  providing explicit relations between the corresponding
fundamental polynomials that will be used afterwards.

In the sequel we will use the convenient abbreviation $T:=Y_{n+1} \cap
K= Y_{n+1}\setminus Y_n$. 

\begin{theorem}\label{T:3.1}
  Let $Y_{n+1}$ be an $\Pi_{n+1}$-poised set and $K$ be a line such that
  $\#(K \cap Y_{n+1})=n+2$.   Then the set $Y_{n}:=Y_{n+1} \setminus K$
  is $\Pi_{n}$-poised.

  Conversely, if $Y_{n}$ is a $\Pi_{n}$-poised set
  with $Y_{n} \cap K=\emptyset$ and  $Y_{n+1} $ is obtained by adding to $Y_{n}$
  $n+2$ distinct nodes on the line $K$, then $Y_{n+1}$ is
  $\Pi_n$-poised. Moreover,
  $$
  \ell_{y,Y_{n+1}} (x) = {k(x) \over k(y)} \ell_{y,Y_{n}}(x), \qquad   y \in Y_{n},
  $$
  and, for $t \in T$,
  \begin{equation}
    \label{eq:(1)}
    \ell_{t,Y_{n+1}}(x) =  \prod_{s \in T \setminus \{t\}} {m(x) -
      m(s)  \over m(t) - m(s)}- k(x) \sum_{y \in Y_{n}}
    {\ell_{y,Y_{n}}(x) \over k(y)} \prod_{s \in T \setminus \{t\}}
    {m(y) - m(s)  \over m(t) - m(s)},
  \end{equation}
  where $m$ is an arbitrary polynomial of first degree such that
  $1,k,m$ form a basis of $\Pi_1$ and $T:= Y_{n+1} \cap K$. 
\end{theorem}

\begin{pf}
  For each $t \in T$, we define the polynomial
  $$
  d_t(x) := \prod_{s \in T \setminus \{t\}} (m(x) - m(s)).
  $$
  Since $d_t(t) \ne 0$, we can find a fundamental polynomial of the form 
  $$
  q_t(x) := {d_t(x) \over d_t(t)}= \prod_{s \in T \setminus \{t\}}
  {m(x) - m(s)  \over m(t) - m(s)}.
  $$
  for each $t \in T$.  So,  $T$ is a $\Pi_{n+1}$-independent set.
  The restriction of the evaluation map 
  $p \mapsto p(T) $ to $\Pi_{n+1}$  is surjective because  $T$ is
  $\Pi_{n+1}$-independent and its kernel $\{p \in \Pi_{n+1}: p(T)
  =0\}$ has dimension  $ \dim \Pi_{n+1} -(n+2)= \dim \Pi_{n}$.
  Clearly $k \Pi_{n} :=\{kp : p \in \Pi_{n}\} $ is contained in its
  kernel and since both spaces have the same dimension they must
  coincide.
  
  Assume that $Y_{n+1}$  is $\Pi_{n+1}$-poised and let $y \in
  Y_{n}:=Y_{n+1}\setminus K$. We have that $\ell_{y,Y_{n+1}}(T)=0$
  which implies $\ell_{y,Y_{n+1}} \in   k\Pi_{n}$.  Therefore, $k$ is a factor
  of $\ell_{y,Y_{n+1}}$ and $k(y) \ell_{y,Y_{n+1}} /k$ is a fundamental
  polynomial in $\Pi_{n}$ of $y$ for the Lagrange interpolation problem
  in $Y_{n}$.  Hence $Y_{n}$ is a $\Pi_{n}$-independent set and since
  $\#Y_{n}=   \# Y_{n+1}  -  (n+2)= \dim \Pi_{n}$, it follows that
  $Y_{n}$ is $\Pi_{n}$-poised, proving the first statement.

  Conversely, if $Y_{n}$ is $\Pi_{n}$-poised, then $k
  \ell_{y,Y_{n}}/k(y)$ is a fundamental polynomial in $\Pi_n$ of $y$ for
  the Lagrange problem in $Y_{n+1}$ for each $y \in Y_{n}$.  Now take
  $t \in T$.   Since $k$ does not vanish on the set $Y_{n}$ and $Y_{n}$
  is $\Pi_{n}$-poised, we can define a polynomial interpolating $q_t/k$
  on $Y_{n}$ 
  $$
  L_{Y_{n}}[q_t/k] = \sum_{y \in Y_{n}} {q_t(y) \over k(y)} \ell_{y,Y_{n}}
  $$
  and deduce that $q_t(x) - k(x)L_{Y_{n}}[q_t/k](x)$
  is a fundamental polynomial for $t$ in $Y_{n+1}$.  Hence $Y_{n+1}$ is
  $\Pi_{n+1}$-independent and $\Pi_{n+1}$-poised since $\#Y_{n+1}=
  \#Y_{n}+ (n+1)= \dim \Pi_{n+1}$.  Finally, (\ref{eq:(1)}) follows from 
  $$
  \ell_{t, Y_{n+1}}(x) =q_t(x) - k(x)L_{Y_{n}}[q_t/k](x)  ={d_t(x) \over
    d_t(t)} - {k(x) \over d_t(t)} \sum_{y \in Y}d_t(y) {\ell_{y,Y}(x)
    \over k(y)}, \qquad t \in T.
  $$
\end{pf}

\begin{theorem}\label{T:3.2}
  Let $Y_n$ be $\Pi_n$-poised, and let $T$ be any set of $n+2$ distinct
  points such that $Y_{n+1} = Y_n \cup T$ is a $\Pi_{n+1}$-poised set.
  Then the $n+2$ functions 
  $$
  h_t:=\ell_{t, Y_{n+1}} \in  \Pi_{n+1}  \cap  I(Y_n), \quad t \in T,
  $$ 
  are  an H-basis of the ideal $I(Y_n)$.   
\end{theorem}

\begin{pf}
  Clearly $\ell_{t, Y_{n+1}} \in \Pi_{n+1}$ vanish on $Y_n$ and so  $
  h_t=\ell_{t, Y_{n+1}} \in \Pi_{n+1}  \cap  I(Y_n)$, $t \in T$.
  Since the Lagrange fundamental polynomials are
  linearly independent, they form a basis of $ \Pi_{n+1}  \cap
  I(Y_n)$.
  From
  the fact that $Y_n$ is $\Pi_n$-poised, it follows that $\Pi = \Pi_n
  \oplus I(Y_n)$ and together with the linear independence of the
  Lagrange fundamental polynomials this yields $\dim \Pi_{n+1} \cap
  I(Y_n) = \dim \Pi_{n+1} - 
  \dim \Pi_n=  n+2$. Now we conclude from Lemma 2.2 that $(h_t: t \in
  T)$ is an H-basis for $I$.
\end{pf}

Theorem~\ref{T:3.1} and Theorem~\ref{T:3.2} can be combined to obtain
an H-basis of the ideal $I(Y)$ of any $\Pi_n$-poised set $Y$ formed by
fundamental polynomials obtained from the Berzolari-Radon construction.

\begin{theorem}\label{T:3.3}
  Let $Y_n$ be a $\Pi_n$-poised set, and $T$ be any set of $n+2$
  points lying on a line $K$ such that $K\cap Y_n=\emptyset$
  and $Y_{n+1} = Y_n \cup T$.   Then the $n+2$ functions 
  $$
  h_t:=\ell_{t, Y_{n+1}} \in  \Pi_{n+1}  \cap  I(Y_n), \quad t \in T,
  $$ 
  are $\RR[k]$-independent functions  and form an H-basis of the ideal
  $I(Y_n)$.
\end{theorem}

\begin{pf}
  By Theorem~\ref{T:3.1}, $Y_{n+1}$ is $\Pi_{n+1}$-poised.  
  The polynomials $\ell_{t, Y_{n+1}} \in \Pi_{n+1}$, $t \in T$, vanish
  on $Y_n$. Let us show now 
  that $ h_t=\ell_{t, Y_{n+1}} \in \Pi_{n+1}  \cap  I(Y_n)$, $t \in
  T$, are $\RR[k]$-independent, that is, if 
  $$
  \sum_{t \in T}  c_t(k) h_t=0,
  $$  
  for some univariate polynomials $c_t$, $t  \in T$, then all these
  polynomials are zero:
  $c_t=0$, $t \in T$.  If we denote by $m_t$ the degree of $c_t$, we can write
  we can write
  $$
  \sum_{t \in T} \sum_{j=0}^{m_t} c_{t,j}k^j h_t=0.
  $$
  After dividing the above equation by an appropriate power of $k$, 
  we may assume that $(c_{t,0}: t \in T) \ne 0$.  Then 
  $$
  \sum_{t \in T} c_{t,0} h_t= - \sum_{t \in T} \sum_{j=1}^{m_t} c_{t,j} k^j h_t =kq
  $$
  for some $q \in \Pi_n$.  Since $k$ does not vanish at any node of
  $Y_n$, we have that $q(Y_n)=0$, that is, $q \in \Pi_n \cap
  I(Y_n)$. But $\Pi_n \cap I(Y_n) =0$ because $Y_n$ is
  $\Pi_n$-poised.  So $q=0$ and from the linear independence of the
  Lagrange fundamental polynomials we deduce that $ c_{t,0}=0$ for
  each $t \in T$.  So, the $\RR[k]$-independence follows.
  The H-basis property follows from Theorem~\ref{T:3.2}.
\end{pf}

\begin{definition}
  A \emph{syzygy} of $P:=(p_0, \ldots, p_{n+1})\in \Pi^{n+2}$ is
  $\Sigma:=(\sigma_0, \ldots, \sigma_{n+1}) \in \Pi^{n+2}$ such that  
  $$
  \sum _{i=0}^{n+1} \sigma_i p_i=0.
  $$
  The set of all syzygies for $P$ will be denoted by $S(P)$ and the
  set of all syzygies of certain maximal degree as $S_k (P) = S(P)
  \cap \Pi_k^{n+2}$, $k \in \NN_0$.
\end{definition}

Inuitively, syzygies describe ambiguities in representing a polynomial
with respect to an ideal. Indeed,
$$
f = \sum _{i=0}^{n+1} f_i p_i = \sum _{i=0}^{n+1} f_i' p_i
$$
holds if and only if $(f_0-f_0',\dots,f_{n+1}-f_{n+1}') \in S(P)$.
The set of $S(P)$ of syzygies of  $P \in \Pi^{n+2}$ forms a
$\Pi$-submodule of $\Pi^{n+2}$.  If $P$ is the basis of an ideal,
$S(P)$ provides important information on the ways of expressing  a
polynomial in the ideal in terms of the basis $P$.

We have seen that for a $\Pi_{n}$-poised set $Y_n$ we can obtain an
H-basis $h_t(x) := \ell_{t,Y_{n+1}}(x)$, $t \in T$, of the ideal
$I(Y_n)$ using the Lagrange fundamental polynomials of a set
$Y_{n+1}=Y_n \cup T$ where  $T$ is a set of $n+2$  points on
a line $K$ such that $K \cap Y_n =\emptyset$.  
We also recall that the function 
$$
d_t(x)= \prod_{s \in T \setminus \{t\}}(m(x) - m(s)) , \quad t \in Y_{n+1},
$$
satisfies
\begin{equation}
  \label{eq:dttneq0}
  d_t (t) \neq 0
\end{equation}
and
can be used to express the Lagrange fundamental polynomials described
in~(\ref{eq:(1)}) 
$$
h_t(x)=\ell_{t,Y_{n+1}}(x) = {d_t(x) \over d_t(t)}- {k(x) \over
  d_t(t)} \sum_{y \in Y_{n}} d_t(y) {\ell_{y,Y_{n}}(x) \over k(y)},
$$
where $m$ is a linear polynomial such that $\Span \{1,k,m\}= \Pi_1$.

First note that by the $\RR[k]$-independence, any nontrivial syzygy
cannot consist of constants or polynomials in $\RR[k]$ only and some
coefficient should include the independent polynomial $m$.
In order to find explicit syzygies, we begin with  $t_i,t_j \in T :=\{t_0,t_1,
\ldots, t_{n+1}\}$ and compute
\begin{eqnarray}
  \nonumber
  \lefteqn{
    d_{t_i}(t_i) (m(x)- m(t_i))h_{t_i}(x)-d_{t_j}(t_j) (m(x)-
    m(t_j))h_{t_j}(x)} \\
  \label{eq:(2)}
  & = & k(x) \sum_{y \in Y_{n}} (d_{t_i}(y) (m(x)-
  m(t_i))-d_{t_j}(y)(m(x)- m(t_j))) {\ell_{y,Y_{n}}(x) \over k(y)}.
\end{eqnarray}
The left hand side of the above equation (\ref{eq:(2)}) belongs to
$I(Y_n)$ and $k$ 
does not vanish at any point of $Y_n$, hence 
$$
\sum_{y \in Y_{n}} (d_{t_i}(y) (m(x)- m(t_i))-d_{t_j}(y)(m(x)-
m(t_j))) {\ell_{y,Y_{n}}(x) \over k(y)}.    
$$
is a polynomial in $\Pi_{n+1} \cap I(Y_n)$.  Since $\{ h_t : t \in T
\}$ is a basis of the vector space $\Pi_{n+1} \cap I(Y_n)$, it follows
that
$$
\sum_{y \in Y_{n}} (d_{t_i}(y) (m(x)- m(t_i))-d_{t_j}(y)(m(x)-
m(t_j))) {\ell_{y,Y_{n}}(x) \over k(y)} = \sum_{s\in T} c_s^{t_i,t_j}
h_s(x),
$$
for coefficients $c_s^{t_i,t_j} \in \RR$.
Since $h_s:= \ell_{s,Y_{n+1}}(x)$ are fundamental polynomials
in $T$, we even have the explicit formula
$$
c_s^{t_i,t_j}= \sum_{y \in Y_{n}} (d_{t_i}(y) (m(s)-
m(t_i))-d_{t_j}(y)(m(s)- m(t_j))) {\ell_{y,Y_{n}}(s) \over k(y)}.
$$
In particular,
\begin{eqnarray*}
  c_{t_i}^{t_i,t_j}
  & = & -(m(t_i)- m(t_j))\sum_{y \in Y_{n}}d_{t_j}(y)
        {\ell_{y,Y_{n}}(t_i) \over k(y)}, \\
  c_{t_j}^{t_i,t_i}
  & = & (m(t_j)-
        m(t_i))  \sum_{y \in Y_{n}} d_{t_i}(y){\ell_{y,Y_{n}}(t_j) \over
        k(y)}.   
\end{eqnarray*}
Then formula (\ref{eq:(2)}) can be written in the form 
$$
d_{t_i}(t_i) (m(x)- m(t_i))h_{t_i}(x)-d_{t_j}(t_j) (m(x)-
m(t_j))h_{t_j}(x)=k(x)\sum_{s\in T} c_s^{t_i,t_j} h_s(x). 
$$
giving rise to syzygies $\Sigma_{t_i,t_j} \in S(H)$ of $H:=(h_t: t \in T)$,
associated to the pairs $t_i, t_j \in T$ whose components
\begin{equation}
  \label{eq:(3)}
  \begin{array}{rcl}
    \sigma_{t_i}^{t_i,t_j}(x) & = & k(x) c_{t_i}^{t_i,t_j} - d_{t_i}(t_i)
    (m(x)-m(t_i)), \\
    \sigma_{t_j}^{t_i,t_j}(x) & = & k(x)c_{t_j}^{t_i,t_j} + d_{t_j}(t_j)
    (m(x)- m(t_j)), \\
    \sigma_s^{t_i,t_j}(x) & = & k(x) c_s^{t_i,t_j},  \qquad s \in T\setminus
    \{t_i, t_j\}, 
  \end{array}
\end{equation}
are polynomials in $\Pi_1$.
 
Note that the syzygies $\Sigma_{t_i,t_j}$, $t_i,t_j \in T$, satisfy 
$$
\Sigma_{t_i,t_j} + \Sigma_{t_j,t_l} = \Sigma_{t_i,t_l}, \quad t_i, t_j, t_l \in T.
$$
In particular, $\Sigma_{t_i,t_i}=0$ and $\Sigma_{t_i, t_j} = - \Sigma_{t_j,t_i}$.  

We now focus on the syzygies
$$
\Sigma_{t_0,t_i}, \quad i=1, \ldots,n+1.
$$
and define, in accordance with \cite{eisenbud05:_geomet_syzyg,Schenck16P},
$$
\Sigma(x)= \left(\sigma_{t_j}^{t_0,t_i}(x)\right)_{i=1, \ldots, n+1, j=0,
  \ldots, n+1} \in \Pi_1^{(n+1)\times(n+2)}. 
$$
as the polynomial matrix whose rows are the components of these syzygies.
We observe that if $x \in K$, then $k(x)=0$ and the $(n+1)
\times(n+1)$ submatrix of $\Sigma(x)$ formed with the $n+1$ last
columns simplifies to a diagonal matrix, whose diagonal entries are 
$$
\sigma_{t_i}^{t_0,t_i}(x)  = d_{t_i}(t_i) (m(x) - m(t_i)).
$$
By definition of $m$ and (\ref{eq:dttneq0}), we have that
$\sigma_{t_i}^{t_0,t_i}(x) \ne 0$ 
for $x \in K \setminus \{t_i\}$, which implies that the rank of
$\Sigma(x)$ over the field of rational functions is $n+1$, in other
words, the syzygies are independent.
Whenever we speak of the rank of a syzygy matrix it has to be
understood in that sense.

The relations $\Sigma(x) H =0$ determine $H$ up to a polynomial factor.   
Let $\Sigma_j(x)$ be the submatrix obtained by removing the column
corresponding to the index $j$, then $\det \Sigma_j(x) \in \Pi_{n+1}$
and there exists $w \in \RR$ such that
$$
h_{t_ j}(x)= (-1)^j \, w \, \det  \Sigma_j(x), \qquad j=0, \ldots, n+1.
$$
In order to determine $w$ we restrict $x$ to the points on the line
$K$, so that $k(x)=0$, and obtain, for $ x \in K$, that
$$
{d_{t_0}(x) \over d_{t_0} (t_0)}=h_{t_ 0}(x)=  w \, \prod_{i=1}^{n+1}
d_{t_i}(t_i) (m(x) - m(t_i))= w \, d_{t_0}(x) \prod_{i=1}^{n+1}
d_{t_i}(t_i),
$$
from which we deduce that 
$$
w = \Big(  \prod_{t \in T} d_t(t)\Big)^{-1} 
$$
and
$$
h_{t_ j}(x)= (-1)^j \Big(  \prod_{t \in T} d_t(t)\Big)^{-1} \det  \Sigma_j(x).
$$
Thus we have shown that the syzygy matrix $\Sigma(x)$ has nonzero
minors $\det \Sigma_j$ for all $j=0, \ldots, n+1$.  Furthermore,
the minors are polynomials of exact degree $n+1$ in $x$.

Let us summarize the results obtained so far for further reference.

\begin{theorem}\label{T:3.4}
  Let $Y_n$ be a $n$-poised set, then there exists for any H-basis  
  $(h_0, \ldots, h_{n+1})$ of $I(Y_n)$ a syzygy matrix
  $\Sigma(x) \in \Pi_1^{(n+1) \times (n+2)} $ of rank $n+1$. The syzygy matrix 
  determines the H-basis up to a constant factor $w \ne 0$
  $$
  h_j(x)= (-1)^j \, w \, \det  \Sigma_j(x), \qquad j=0, \ldots, n+1,
  $$
  where $\Sigma_j(x)$ denotes the submatrix obtained from $\Sigma(x)$ by
  removing the column corresponding to the index $j$.
\end{theorem}

\begin{pf}
  We have constructed $\Sigma(x)$ for a particular H-basis.  Since all
  H-bases with $n+2$ elements of $I(Y_n)$ are bases of $\Pi_{n+1} \cap I(Y_n)$, they are
  related by a nonsingular matrix in $B\in \RR^{(n+2) \times (n+2)}$.
  Multiplying $\Sigma(x)$ from the right with $B$, we obtain the
  syzygy matrix $\Sigma(x) B  \in \Pi_1^{(n+1) \times (n+2)}$ of rank
  $n+1$ for an arbitrary H-basis $B^{-1} H$.   
\end{pf}


\begin{theorem}\label{T:3.5}
  Let $Y_n$ be a $n$-poised set.  If $\Sigma(x), \Sigma'(x) \in
  \Pi_1^{(n+1) \times (n+2)}$ are syzygy matrices of rank $n+1$ for
  two H-bases of $I(Y_n)$, then there exist nonsingular
  scalar
  matrices $A$ and $B$, such that 
  $$
  \Sigma'(x) = A \Sigma(x) B,
  $$
  i.e., the linear syzygy matrix of rank $n+1$ for $I(Y_{n})$ is unique
  up to equivalence.
\end{theorem}

\begin{pf}
  Let us denote by  $S_1(H) =S(H) \cap \Pi_1^{n+2}$ the space of
  linear syzygies of a given H-basis of $I(Y_n)$. By Theorem~\ref{T:3.4},
  there exists a linear syzygy matrix $\Sigma(x)$ with respect to the
  particular H-basis
  $$
  H= ( h_j(x):=E_{Y_n}[x_1^{n+1-j}x_2^j]: j =0, \ldots, n+1),
  $$
  where $E_{Y_n}[f]$ is the interpolation error to $f$ on $Y_n$ from $\Pi_n$.
  We will now show that any syzygy in $S_1(H)$  can be written
  as a linear combination of the rows of $\Sigma(x)$.  Since
  $$
  x_2 E_{Y_n}[x_1^{n+1-i}x_2^i](x) - x_1
  E_{Y_n}[x_1^{n-i}x_2^{i+1}](x) = \sum_{j=0}^{n+1}c_{ij}
  E_{Y_n}[x_1^{n+1-j}x_2^j](x), 
  $$
  the corresponding syzygy matrix for $H$ can be written as
  $$
  \Sigma(x)=\pmatrix{
    c_{00} -x_2 & c_{01} + x_1 & c_{02} & \dots  & c_{0, n+1} \cr
    c_{10} & c_{11} -x_2 & c_{12} + x_1  & \dots & c_{1, n+1} \cr
    \vdots & \ddots & \ddots &\ddots & \vdots \cr
    c_{n0} & \dots  & c_{n,n-1} &c_{n, n}-x_2  & c_{n, n+1} + x_1 \cr
  }
  $$
  Now assume that $Z(x)$ is another syzygy in $S_1(H)$.  We subtract
  a proper multiple of the first row of $\Sigma(x)$ from $Z(x)$ in
  such a way that the second component does not depend on $x_1$, then
  use the second row to eliminate $x_1$ from the third component, and
  so on.  The resulting syzygy $Z(x)-a^T\Sigma(x)$ is of the form 
  $$
  (\sigma_0(x_1,x_2), \sigma_1(x_2), \ldots, \sigma_{n+1}(x_2)).
  $$
  But then 
  $$
  0= \sum_{j=0}^{n+1} \sigma_j(x) h_j(x) =\sigma_0(x_1,x_2) E_{Y_n}
  [x_1^{n+1}](x) + \sum_{j=1}^{n+1} \sigma_j(x_2)
  E_{Y_n}[x_1^{n+1-j}x_2^j](x) 
  $$
  implies that $\sigma_0 =0$ because it is the coefficient of the only
  appearance of the power $x_1^{n+1}$. We can inductively apply
  the same reasoning to the powers $x_1^{n+1-j} x_2^j$ to conclude
  that $\sigma_j=0$, $j=1, \ldots, n$.  
  Therefore $Z(x) = a^T \Sigma(x)$, that is, $Z(x)$ is a linear combination
  of the rows of $\Sigma(x)$, hence $\dim S_1(H)= n+1$. 
 
  So we know that two linear syzygy matrices of rank $n+1$ for the
  same H-basis are related by left multiplication by a nonsingular
  matrix.  On the other hand, the proof of Theorem~\ref{T:3.4} tell us that
  changes of the H-bases correspond to right multiplication by a
  nonsingular matrix.   
\end{pf}

\section{Maximal lines}
Maximal lines, or, more generally, maximal hyperplanes, as introduced
in \cite{boor07:_multiv},  are at the heart of the Gasca-Maeztu
conjecture. A $\Pi_n$-poised set $Y_n$ is said to  
contain a \emph{maximal line} if there exists a line $K$ such that 
$\# ( Y_n \cap K ) = n+1$. As pointed out first by H.~Schenck, the
following result that closely connects maximal lines to the syzygy
matrix $\Sigma(x)$, can be seen as a special case of the
Hilbert--Burch theorem, cf. \cite{eisenbud05:_geomet_syzyg}.
We restate this fact here and give a more direct, affine and
elementary proof of it. 

\begin{proposition}\label{P:4.1}
  Suppose that $Y_n$ is a $\Pi_n$-poised set and 
  there exist linearly independent $h_0,\dots,h_{n+1} \in \Pi_{n+1}
  \cap I(Y_n)$ with a  syzygy matrix $\Sigma \in \Pi_1^{(n+1) \times
    (n+2)}$  of rank $n+1$ with one column  of the form $k(x) v$, $v\in\RR^{n+1}$, for some nonconstant $k \in \Pi_1$. Then $K=V(k)$ is a maximal line for $Y_{n}$,
  i.e., $\# ( K \cap Y_n ) = n+1$. 
\end{proposition}

\begin{pf}
  Since the polynomials $h_0, \ldots, h_{n+1}$ generate $I(Y_n)$ and
  $\dim (\Pi_{n+1} \cap I(Y_n)) =n+2$, it follows from Lemma \ref{L:2.2} that
  these polynomials form an $H$-basis of $I(Y_n)$. 
  After renumbering the ideal basis, we can assume that the last column of $\Sigma(x)$ is of the form $k(x) v$ for some vector $v$.  By Theorem \ref{T:3.4},
  $h_j =w (-1)^j\det \Sigma_j \in \Pi_{n+1}$, $j=0,\dots,n+1$, up to a
  nonzero constant factor $w$.  It follows that $h_j (x) := k(x) g_j (x)$, $g_j
  \in \Pi_n$, $j=0, \ldots, n$, 
  all belong to $I( Y_{n} )$ and that $g_j \in I ( Y_n \setminus K )$,
  $j=0,\dots,n$, are $n+1$ linearly independent polynomials of degree
  $n$ in $I ( Y_n \setminus K )$.  If $q \in I(Y_n \setminus K) \cap
  \Pi_{n-1}$, then $qk \in \Pi_n$ vanishes on $Y_n$.  Since $Y_n$ is
  $\Pi_n$-poised, it follows that $q=0$.  So we have $I(Y_n \setminus
  K) \cap \Pi_{n-1}=0$.  By Lemma~\ref{L:2.2}, $(g_j: j=0, \ldots,n)$
  is an $H$-basis of $I ( Y_n \setminus K )$.  Since $\Pi = \Pi_{n-1}
  \oplus I ( Y_n \setminus K )$, we deduce that $Y_n \setminus K$ is
  $\Pi_{n-1}$-poised and  
  $$
  \# ( Y_n \cap K ) = \dim \Pi_n - \# ( Y_n \setminus K ) 
  = \dim \Pi_n - \dim \Pi_{n-1} = n+1
  $$
  which means that $K$ is indeed a maximal line.  
\end{pf}

\noindent
Since one can apply row transformations to the syzygy matrix to obtain another syzygy matrix for the same basis, the following statement is equivalent to
Proposition~\ref{P:4.1}. 

\begin{corollary}\label{C:4.2}
  Suppose that $Y_n$ is a $\Pi_n$-poised set and 
  there exist linearly independent $h_0,\dots,h_{n+1} \in \Pi_{n+1}
  \cap I(Y_n)$ with a  syzygy matrix $\Sigma \in \Pi_1^{(n+1) \times
    (n+2)}$  of rank $n+1$ whose $j$-th column is of the form $k(x) \,
  e_i$ for some linear function $k$ and some $i \in \{ 1,\dots,n+1 \}$
  and $j \in \{ 0,\dots,n+1 \}$. Then $K$ is a maximal line for
  $Y_{n}$. 
\end{corollary}

The results from the preceding section also allow us to give a converse 
statement of Proposition~\ref{P:4.1} that says that any maximal line
can be found in a proper syzygy matrix.

\begin{proposition}\label{P:4.1c}
  If a $\Pi_n$-poised set $Y_n \subset \RR^2$ contains a maximal line $K$, then
  there exists a basis  $h_0,\dots,h_{n+1}$ of $\Pi_{n+1} \cap I(Y_n)$
  with a  syzygy  
  matrix $\Sigma \in \Pi_1^{(n+1) \times (n+2)}$ of rank $n+1$ whose
  last column is of the form $k(x) \, e_{n+1}$.
\end{proposition}

\begin{pf}
  Let $Y_{n-1}$ be any $(n-1)$-poised 
  set in $\RR^2$ and $K$ be a line with $K \cap Y_{n-1} = \emptyset$. Choose 
  $n+1$ points $t_0,\dots,t_n$ on $K$ and let $Y_n$ be the union of $Y_{n-1}$ and 
  these $n+1$ points. Now set, with the notation from the preceding section,
  $$
  g_j (x) := k(x) \, h_{t_j} (x), \quad j=0,\dots,n, \qquad
  g_{n+1} (x) :=- ( m(x) - m(t_0) ) h_{t_0} (x).
  $$
  These polynomials are linearly independent. Indeed, evaluating
  $$
  \sum_{j=0}^{n+1} c_j \, g_j (x) = 0
  $$
  along $K$ yields $c_{n+1} = 0$ and then
  $$
  0 = k(x) \, \sum_{j=0}^{n} c_j \, h_{t_j} (x),
  $$
  which implies that $c_0 = \cdots = c_n = 0$. Since, in addition,
  $g_j (Y_n) = 0$, it follows from Lemma~\ref{L:2.2} that 
  $g_0,\dots,g_{n+1}$ are an H-basis of $I(Y_n)$.  The syzygy matrix for this 
  basis takes the form
  $$
  \left( 
    \matrix{ 
      \sigma^{t_0,t_1}_{t_0} (x) & \sigma^{t_0,t_1}_{t_1} (x) & \dots 
      & \sigma^{t_0,t_1}_{t_n} (x) & 0 \cr
      \vdots & \vdots & \ddots & \vdots & \vdots \cr
      \sigma^{t_0,t_n}_{t_0} (x) & \sigma^{t_0,t_n}_{t_1} (x) &\dots 
      & \sigma^{t_0,t_n}_{t_n} (x) & 0 \cr
      m(x) - m(t_0) & 0 & \dots & 0 & k(x)}
  \right),
  $$
  with $\sigma_{t_j}^{t_0,t_i}$ given by (\ref{eq:(3)}). The upper left
  part of this matrix is the syzygy matrix for $Y_n \setminus K$.
  The last column now consists 
  of the linear polynomial $k$ marking the maximal line.  
\end{pf}

Combining Proposition~\ref{P:4.1}, Proposition~\ref{P:4.1c} and
Theorem~\ref{T:3.5}, we can now even give a 
characterization of maximal lines in terms of syzygies.

\begin{theorem}\label{T:4.3}
  A $\Pi_n$-poised set $Y_n \subset \RR^2$ contains a maximal line if
  and only if
  there exists a basis  $h_0,\dots,h_{n+1}$ of $\Pi_{n+1} \cap I(Y_n)$
  with a  syzygy  
  matrix $\Sigma \in \Pi_1^{(n+1) \times (n+2)}$ of rank $n+1$ whose
  $j$-th column is of the form 
  $k(x) \, e_i$ for some  polynomial $k$ of degree 1 and some $i \in \{ 1,\dots,n+1 \}$ 
  and $j \in \{ 0,\dots,n+1 \}$. The  polynomial $k$ determines the maximal line 
  $K = V(k)$.
\end{theorem}

\section{Syzygy matrices of $GC_n$ sets}

We recall that a $GC_n$ set, introduced by Chung and Yao
\cite{ChungYao77}, is a $\Pi_n$-poised set whose Lagrange fundamental
polynomials are products of linear factors.  Gasca and Maeztu
conjectured in  \cite{GascaMaeztu82} that any $GC_n$ set contains a
maximal line.  
In this section we consider some special $GC_n$ sets and their syzygy
matrices, namely
the two most important and best investigated
\cite{deBoor97,CarnicerGascaSauer09}  examples of $GC_n$ sets: natural
lattices and (generalized) principal lattices.
We start with the \emph{natural lattice} $Y_n$ corresponding to the
intersections of $n+2$ lines $K_0, \ldots, K_{n+1}$ in general
position, such that
$K_i  \cap K_j$, $0 \le i < j \le n+1$, form a set of ${n+2 \choose
  2}$  points $x_{ij}$ in $\RR^2$.  The
Lagrange fundamental polynomial $\ell_{ij}$ corresponding to $x_{ij}$
is then of the form 
$$
\ell_{ij} (x) = \prod_{r\in \{0, \ldots, n+1\} \setminus\{ i,j\}}
{k_r(x) \over k_r(x_{ij})}. 
$$
Then we can construct an H-basis of $I(Y_n)$ basis  using
Theorem~\ref{T:3.2}.  Let $K_{n+2}$ be a line intersecting all
previous lines at $x_{i,n+2}$, then 
$$
\ell_{i,n+2}(x) =  \prod_{r\in \{0, \ldots, n+1\} \setminus\{ i\}}
{k_r(x)  \over k_r(x_{i,n+2})}. 
$$
So,
$$
h_i(x) = \ell_{i,n+2}(x) \, \prod_{r \ne i}  k_r(x_{i,n+2}) =  \prod_{r\in
  \{0, \ldots, n+1\} \setminus\{ i\}} k_r(x) , \qquad i=0, \ldots,
n+1, 
$$
form an H-basis of $I(Y_n)$. For this H-basis we obtain the following
syzygy matrix:
$$
\Sigma(x) = \pmatrix{
k_0(x) &- k_1(x)  & 0 & \dots & 0 \cr
k_0(x) & 0 &- k_2(x) & \ddots &\vdots \cr
\vdots & \vdots & \ddots & \ddots & 0 \cr
k_{0}(x) & 0 & \dots & 0 & -k_n(x) \cr}.
$$
Corollary~\ref{C:4.2} proves that
$K_j$, $j= 0 \ldots, n+2$, are maximal lines. 

Now let us consider the case of a generalized principal lattice,
introduced in \cite{CarnicerGasca05, CarnicerGasca06} and further
analyzed in \cite{CarnicerGodes06}.  In 
this case, we have $3n$ lines $K_{i,j}$, $i=0, \ldots, n$, $j=0,1,2$,
where 
$$
K_{\beta_0,0} \cap K_{\beta_1, 1} \cap K_{\beta_2,2} = \{x_\beta\},
\quad \vert \beta \vert =n. 
$$
are requested to be distinct.  Then $Y_n= \{x_\beta : \vert \beta
\vert = n\}$ is an $\Pi_n$-poised set and $GC_n$ because the
corresponding Lagrange polynomials are of the form 
$$
\ell_\beta (x) = \prod_{\gamma_0< \beta_0} {k_{\gamma_0,0} (x) \over
  k_{\gamma_0,0} (x_\beta)} \prod_{\gamma_1< \beta_1} {k_{\gamma_1,1}
  (x) \over k_{\gamma_1,1} (x_\beta)}\prod_{\gamma_2< \beta_2}
{k_{\gamma_2,2} (x) \over k_{\gamma_2,2} (x_\beta)}. 
$$
One H-basis of $I(Y_n)$ is
\begin{equation}
  \label{eq:PrincLatHBasis}
  h_j(x) = \prod_{\gamma_1< j} k_{\gamma_1,1} (x)\prod_{\gamma_2< n+1-j}
  k_{\gamma_2,2} (x) , \qquad j=0, \ldots, n+1,   
\end{equation}
whose syzygy matrix is of the form
$$
\Sigma(x) = \pmatrix{
k_{0,1}(x) &-k_{n,2}(x)  & 0 & \dots & 0 \cr
0  &k_{1,1} (x)& -k_{n-1,2}(x) & \ddots &\vdots \cr
\vdots & \ddots & \ddots & \ddots & 0 \cr
0& \dots & 0 & k_{n,1}(x) &  -k_{0,2}(x) \cr}.
$$
Thus $K_{0,1}$ and $K_{0,2}$ turn out to be maximal lines by
inspecting first and last column of the syzygy matrix. Since we left out any
factor with $k_{j,0}$ in the H--basis (\ref{eq:PrincLatHBasis}), the
line $K_{0,0}$ is not detected by the syzygy matrix. However, there
must exist a basis transform that maps the polynomials from
(\ref{eq:PrincLatHBasis}) to 
\begin{equation}
  \label{eq:PrincLatHBasis2}
  h_j'(x) = \prod_{\gamma_0< j} k_{\gamma_0,0} (x) \prod_{\gamma_2< n+1-j}
  k_{\gamma_2,2} (x) , \qquad j=0, \ldots, n+1,   
\end{equation}
which implies that a column transform of $\Sigma (x)$ then gives a column
consisting only of a multiple of $k_{0,0}$.

The particular structure of the fundamental polynomials of $CG_n$
sets, namely the rare property that the can be factorized into linear
polynomials, 
and the preceding examples suggest the following construction of a
factorizable
H-basis for $I (Y_n)$ which is originally due to Schenck in an
unpublished manuscript and for which we can now
give a more direct and elementary exposition. To this end, we choose a
point $z \in \RR^2$ such that $\ell_{y,Y_n} (z) \neq 0$, $y \in Y_n$, and
 define the sets
$$
\cL (y) := \{ m \in \Pi : m | \ell_{y,Y_n}, m(z) = 1 \}, \quad y \in
Y_n,
\qquad \cL := \bigcup_{y \in Y_n} \cL (y).
$$
The purpose of the point $z$ is only to uniquely normalize the
polynomials in $\cL$. Moreover, let $Y_{n-1}$ be any $\Pi_{n-1}$-poised subset
of $Y_n$ and set $T := Y_n \setminus Y_{n-1}$.

\begin{proposition}\label{P:HalBasis}
  The set
  \begin{equation}
    \label{eq:HalBasis}
    \bigcup_{t \in T}
    \left\{ m(x) \, \ell_{t,Y_n} (x) : m \in \cL, m(t) = 0 \right\} 
  \end{equation}
  contains an H-basis for $I( Y_n )$ whose elements are product of first degree polynomials.
\end{proposition}

\begin{pf}
  Since,
  there always exist at least two
  nonparallel lines belonging to $\cL$ that vanish on $t$, we find
  that
  $$
  \Span \left\{ m(x) \, \ell_{t,Y_n} (x) : m \in \cL, m(t) = 0
  \right\} = \Span \left\{ x_1 - t_1, x_2 - t_2 \right\}
  \ell_{t,Y_n} (x).
  $$
  Since, on the other hand,
  $$
  \Span \{ \ell_{t,Y_n} (x) : t \in T \} = \Span \{ 
 E_{Y_{n-1}} [x^\alpha] : |\alpha| = n \},
  $$
  it follows that the space
  $$
  \sum_{t \in T} \Span \left\{ m(x) \, \ell_{t,Y_n} (x) : m
    \in \cL, m(t) = 0 \right\} \subset I(Y_n)
  $$
  contains all polynomials of the form $E_{Y_{n}} [x^\alpha]$,
  $|\alpha| = n+1$.
  Indeed, write, for some $\alpha$ with $|\alpha| = n$
  $$
E_{Y_{n-1}} [x^\alpha] = \sum_{t \in T} c_t \, \ell_{t,Y_n} (x),
  $$
  and use $\epsilon_1 := (1,0)$ and $\epsilon_2 := (0,1)$ for the unit
  multiindices in $\NN_0^2$. It then follows for $j=1,2$, that
  $$
    x_j \, E_{Y_{n-1}} [x^\alpha] - \sum_{t
    \in T} c_t \, t_j \, \ell_{t,Y_n} (x) = \sum_{t
    \in T} c_t \, (x_j - t_j ) \ell_{t,Y_n} (x),
  $$
  is a  polynomial of degree $n+1$ which vanishes on $Y_n$.   This polynomial coincides with 
  $ E_{Y_n} [x^{\alpha+\epsilon_j}]$ because 
  $x^{\alpha+\epsilon_j}-x_j E_{Y_{n-1}} [x^\alpha]  \in \Pi_n$.  So, we have
  $$
  E_{Y_n}
  [x^{\alpha+\epsilon_j}] = \sum_{t
    \in T} c_t \, (x_j - t_j ) \ell_{t,Y_n} (x)
  $$
  and therefore the set from (\ref{eq:HalBasis}) spans $\Pi_{n+1} \cap
  I(Y_n)$. By 
  Corollary~\ref{C:HBasisContained} the set is 
  an H-basis as claimed. Since all candidates are factorizable by
  construction, so is the resulting H-basis.  \end{pf}

\begin{remark}
  Proposition~\ref{P:HalBasis} can also be found in \cite{Schenck16P}
  where it is in turn attributed as implicitly given already in
  \cite{SauerXu95}.
\end{remark}

\section*{Acknowledgement}
We a very grateful to Hal Schenck for the discussions on a draft of
his paper that helped us to write this paper and to Carl de Boor for a
lot of valuable remarks and comments on the first draft that greatly
improved the presentation.


\providecommand{\bysame}{\leavevmode\hbox to3em{\hrulefill}\thinspace}
\providecommand{\MR}{\relax\ifhmode\unskip\space\fi MR }
\providecommand{\MRhref}[2]{%
  \href{http://www.ams.org/mathscinet-getitem?mr=#1}{#2}
}
\providecommand{\href}[2]{#2}

\end{document}